# A short-term planning model for high-speed train assignment and maintenance scheduling

Boliang Lin[*]

School of Traffic and Transportation, Beijing Jiaotong University, No. 3 Shangyuancun, Haidian District, Beijing 100044, China

**Abstract**

This paper considers the simultaneous high-speed train assignment and maintenance scheduling problem. For the maintenance scheduling module, our focus is the second-level maintenance, which is carried out once a month on average. We propose a binary non-linear programming model to mathematically describe the problem. The objective aims to minimize the mileage losses for all high-speed trains and the constraints cover various operational requirements and capacity restrictions. To model the difficult operational requirements, a novel cumulative mileage update function is developed; meanwhile, to describe the depot maintenance capacity restriction, we employ a state function that is able to identify whether a train is in the maintenance state or in the operation state.

**Keywords:** High-speed railway; rolling stock assignment; maintenance scheduling; binary non-linear programming

## 1. Introduction

In recent years, the high-speed railway in China develops quickly, and the number of high-speed trains that put into operation has been increasing over time. By the end of year 2016, the total operation mileage of China railway was 124, 000 kilometers, and there were 2,586 (standard) high-speed trains serving for the high-speed railway passenger transportation. Therefore, in order to further improve the service for high-speed railway transportation organization, the problem that how to enhance the high-speed train management including train assignment and maintenance scheduling is becoming a focus of attention for the railway transportation department. However, the problem is difficult to be solved effectively. High-speed trains have a few special characteristics compared with the general rolling stock. Firstly, the procurement cost of a High-speed train is significantly high, for example a normative "CRH" train is worth of 200,000,000 RMB approximately. Secondly, the inspection and repair system of high-speed train in China is extremely complex, which includes a lot of different maintenance packets (items). Thirdly, the mileage of different routes ranges a lot, and there is a strong relevance between high-speed train assignment and maintenance scheduling, and so on. Because of this, it is necessary to take the train assignment and maintenance scheduling into consideration simultaneously while developing the high-speed train circulation plan, which aims at improving the quality of the plan. In this way, it not only increases the high-speed train operation efficiency and reduces the quantities of high-speed train as well as the procurement cost, but also reduces the maintenance times and the maintenance cost during a fixed period.

---

[*] E-mail: bllin@bjtu.edu.cn

This paper considers the simultaneous high-speed train assignment and maintenance scheduling problem. For the maintenance scheduling module, our focus is the second-level maintenance, which is carried out once a month on average. We propose a binary non-linear programming model to mathematically describe the problem. The objective aims to minimize the mileage losses for all high-speed trains and the constraints cover various operational requirements and capacity restrictions. To model the difficult operational requirements, a novel cumulative mileage update function is developed; meanwhile, to describe the depot maintenance capacity restriction, we employ a state function that is able to identify whether a train is in the maintenance state or in the operation state.

The reminder of this paper is structured as follows: Section 2 presents a comprehensive literature review on related studies. In Section 3, we introduce the problem of high-speed train assignment and maintenance scheduling. In Section 4, we analyze the optimization objective and constraints, and a binary non-linear programming model is proposed. Finally, conclusions are drawn and future research directions are discussed in Section 5.

## 2. Literature Review

With the development of high-speed railway around the world, many experts and scholars have researched the high-speed train assignment and maintenance scheduling problem from different aspects, and proposed different models and solution methods for the high-speed train circulation. In the existing literature, Kroon et al. [1～4] researched the high-speed train assignment and connection problem extensively and profoundly, these researches mainly aimed at minimizing the quantities of high-speed trains and enhancing the robustness of the high-speed train circulation plan which was based on the given train arcs and timetable, and considered the constraints including the types of high-speed trains and the maintenance conditions, etc. Then a few diverse mathematical integer programming models were proposed according to the different concrete issues, and some algorithms such as branch and bound method were designed to solve the optimization models. Cordeau et al. [5] proposed a multi commodity network flow-based model for assigning locomotives and cars to trains in the context of passenger transportation, and adopted a branch-and-bound method to solve the problem. Lingaya et al. [6] researched locomotives and cars assignment to a set of scheduled trains for the passenger railways, and then described a modeling and solution methodology for a car assignment problem. Noori and Ghannadpou [7] studied the locomotive assignment problem which is modeled using vehicle routing and scheduling problem in their research, and a two-phase approach based on a hybrid genetic algorithm was used to solve the problem. Maróti and Kroon [8, 9] focused on the regular preventive maintenance of train unites at NS Reizigers, and presented two integer programming models for solving the maintenance routing problem, one is the interchange model, and the other is the transition model, which could solve the maintenance issue in the forth-coming one to three days for the train units. Giacco [10] focused on the high-speed train connection problem on a transportation service network, and then aimed at minimizing the quantities of high-speed train and constructed a mixed integer linear programming model for optimizing the short-period maintenance scheduling of high-speed trains, and the model took the transportation tasks, running without taking passengers, short-period maintenance items into consideration. Many experts and scholars have researched the same problem in other relevant industry such as the fleet assignment problem (FAP), etc. Lai et al. [11] developed an exact optimization model to improve the efficiency in rolling stock usage with consideration of all necessary regulations and practical constraints. A hybrid heuristic process was also developed to improve solution quality and efficiency. Subsequently, they extended their work into high-speed railway system, where an exact optimization model and a heuristic method were proposed to automate the train-set rostering planning process for Taiwan High-Speed Rail [12]. Moudani and Mora-Camino [13] proposed a dynamic approach for the problems of assigning planes to fights and of feet maintenance operations scheduling in their research. Sherali et al. [14] presented a method to integrate the FAP with schedule design, aircraft maintenance routing, and crew scheduling, and presented a randomized search procedures. Deris et al. [15] researched the problem of ship maintenance scheduling and modelled it as a constraint satisfaction problem (CSP) in their paper, a genetic algorithm (GA) was adopted to solve the problem. Go et al. [16] researched the problem of operation and maintenance scheduling for a containership, and developed a mixed integer programming model for the problem, based on which a heuristic algorithm was presented.

There are also many Chinese specialists and scholars who have researched the relative problem according to the actual situation in China. Wang et al. [17] analyzed three status of High-speed train, including undertaking route, being in maintenance and waiting for maintenance, and then a connection network composed of undertaking route and conducting maintenance is designed, based on which the author proposed an optimization model aiming at maximizing the accumulating mileage before conducting the corresponding maintenance. Wang et al. [18] researched the integrated optimization method for operation and maintenance planning of high-speed train, which aimed at reducing the quantities and the maintenance cost of high-speed train, and designed a max-min ant colony algorithm to solve the problem. Lin et al. [19] described three different operation modes of high-speed trains, and focused on the characteristics of the Grade 2 Maintenance on uncertain railroad region and certain railroad regions. Two optimization models of high-speed train maintenance plan were constructed under these two different operation modes of the trains.

From the literatures mentioned above, we could get a conclusion that the related issues about high-speed train circulation problem are researched extensively and profoundly, and a few research achievements are got, such as optimization models and solution algorithms. However, the researches and the achievements are focused on some specific issues on the whole. What's more, the condition of high-speed railway in China is not the same with abroad, and the corresponding research achievements abroad could not be used to serve the transportation organization of high-speed railway in China directly, and those in China mostly stay at a level of theoretical research. All the research gaps motivate this paper.

## 3. Problem Statement

The issue discussed in this paper is the high-speed train assignment and maintenance scheduling, and it refers to a few key elements including high-speed trains, routes, and maintenance packets, etc. Therefore, the main tasks of high-speed train assignment and maintenance scheduling are to assign a well-conditioned high-speed train to each route every day, and to arrange the maintenance work of high-speed train, of which the accumulated mileage or time of the corresponding maintenance item after the last same maintenance is to meet the maintenance period. In this paper, we define the transportation tasks (train routes) as the ordered trains' circulation which are undertook by the same high-speed train from the departure to the arrival at the high-speed train depot. That is to say, the departure depot and the arrival depot is the same depot which is the attachment depot of the high-speed train. In China, there is a certain time reserved for infrastructure inspection of high-sped railway, which is general four to six hours, and the train is forbidden during this period. So the high-speed train regularly stays at the attachment depot or the other depot at night. Therefore, according to the departure time and the arrival time of each route, we divide the route into one-day route and multi-day route in this paper. In our opinion, if the departure time and the arrival time of a route is on the same day (00:00～24:00), we call the route an one-day route. And if the arrival time of a route is on the next day relative to the departure time, we call the route a two-day route, and the others can be called as multi-day route by that analogy.

According to the definition of route, we give a diagrammatic sketch of route which is shown in Fig. 1 (see also Li et al. [20]). In the diagrammatic sketch, there are three stations named station A, station B and station C respectively, and the high-speed train depot is at the same place with station B. The high-speed trains depart from station B, and undertake the trains from station B to station A and from station B to station C respectively, and come back to station B with undertaking the corresponding trains, and the high-speed trains may go to the high-speed train depot for fixed maintenance or staying. So the route composed of train T1 and train T2 is a one-day route, and the route composed of train T3 and train T4 is a two-day route.

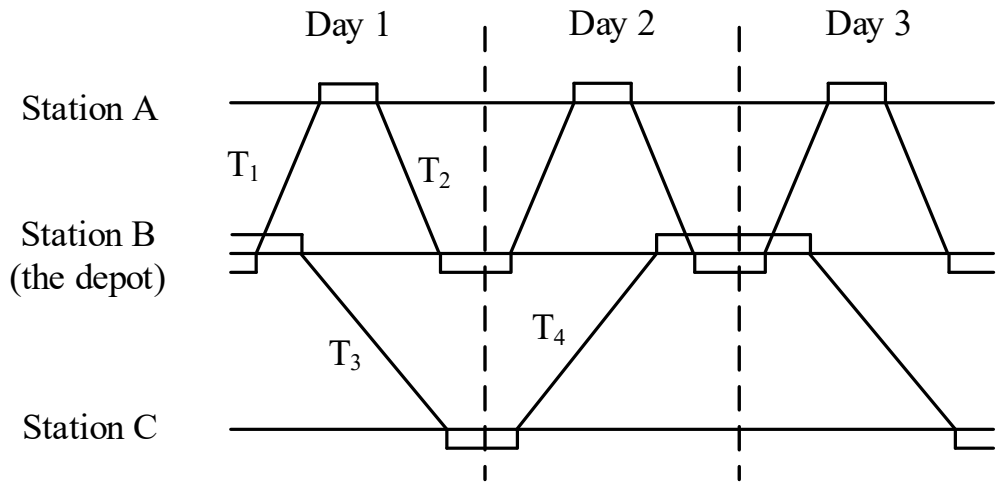

**Fig. 1.** A diagrammatic sketch of route.

It is relatively complex for the inspection and repair system of high-speed train in China and the maintenance contents are divided into five levels. The first-level maintenance is generally called the routine maintenance, the second-level maintenance is generally called the special maintenance, and the third- to fifth-level maintenance are generally called the high-level maintenance uniformly. While the accumulated mileage is to meet 4,000 kilometers or the accumulated time is to meet 48 hours after the last routine maintenance, the high-speed train have to go back to the depot for the routine maintenance again. Compared with other two kinds of maintenance, the maintenance period of the routine maintenance is shorter and it occurs generally at night, so it is of high maintenance frequency. On the contrary, the high-level maintenance has a longer maintenance period and the time spent on the maintenance is comparative longer as well. The special maintenance has a few maintenance packets, such as I2 maintenance, M1 maintenance, flaw detection of hollow axle, traction engine greasing, and so on. Besides, each maintenance item of the special maintenance ranges a lot in the aspects of maintenance period (including mileage and time cycle) and maintenance service time. Therefore, the problem of second-level maintenance scheduling is much more complex than any other, and its quality has a deep effect on the operation efficiency and the maintenance cost of high-speed trains. Because of this, we focus on taking the second-level maintenance into consideration in this paper, and research the optimization method for train assignment and maintenance scheduling. In Fig. 2 (see also Li et al. [20]) we give an example of high-speed train assignment and maintenance schedule with consideration of the second-level maintenance.

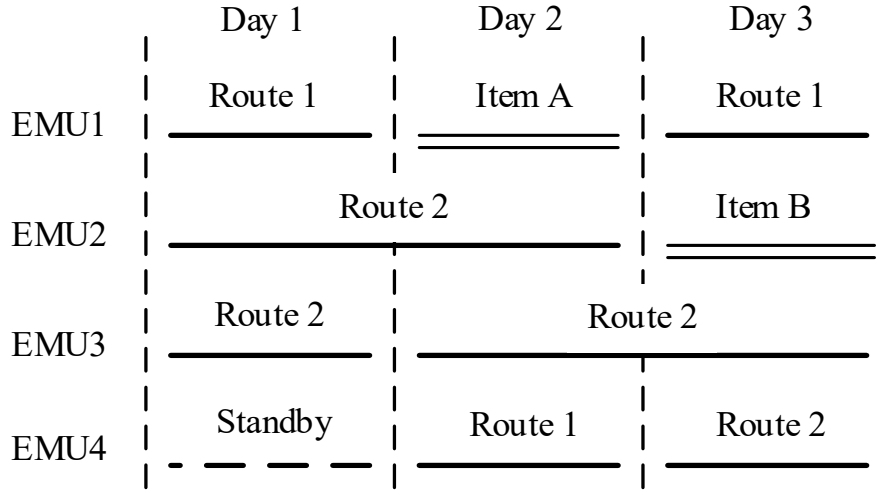

**Fig. 2.** An example of high-speed train assignment and maintenance schedule.

For the high-speed train assignment and maintenance schedule shown in Fig.2. The high-speed train EMU1 undertakes route 1 on the first day, and conducts maintenance item A on the second day, and undertakes route 1 again on the third day. The high-speed train EMU2 undertakes route 2 on the first two days, and conducts maintenance item B on the third day. The high-speed train EMU3 undertakes the route 2 during the period. The high-speed train EMU 4 is being in standby state on the first day, and undertakes the route 1 and route 2 respectively on the second day and the third day.

## 4. Mathematical Model

In this section, we first make some basic assumptions to facilitate the model formulation, followed by detailed descriptions of the notations used in this paper. After that, a binary non-linear programming model is proposed for the high-speed train assignment and second-level maintenance planning problem.

**4.1 Basic assumptions**

To facilitate the model formulation, we make several assumptions throughout this paper which are explained as follows:

- **Assumption 1:** We assume that within the planning horizon a high-speed train carries out the second-level maintenance at most once. This assumption is reasonable because as per the current practice in China high-speed railway system, the planning period is usually one week, which means it is far shorter than the cycle of second-level maintenance.
- **Assumption 2:** We assume the cycles of second-level maintenances are identical, regardless of involved maintenance packets. Although in practice the cycle of second-level maintenance is packet-dependent, given Assumption 1, we could replace the maintenance cycle with an accurate value in the next planning period.
- **Assumption 3:** To the best of our knowledge, there exits no train route that lasts for three or more than three days. Therefore, we assume a train route can be decomposed into at most two sub-routes.

## 4.2 Notations

The notations used in this paper are listed in Table 1.

**Table 1** Notations used in this paper.

| Symbol | Description |
|---|---|
| **Sets** | |
| $E$ | Set of high-speed trains. |
| $P$ | Set of second-level maintenance packets (items). |
| $R$ | Set of transportation tasks (train routes). |
| **Indices** | |
| $m$ | Index of trains. |
| $p$ | Index of second-level maintenance packets, $p \in P$. |
| $r$, $s$ | Index of transportation tasks, $r, s \in R$. Note that when a train is in a spare state, $r, s = 0$. |
| $t$ | Index of days in the considered planning horizon. |
| **Parameters** | |
| $T$ | Planning horizon. |
| $b_{mp}$ | $b_{mp} = 1$ if the next maintenance packet of train $m$ is $p$; $b_{mp} = 0$ otherwise. |
| $a_{mr}$ | $a_{mr} = 1$ if train $m$ can execute transportation task $r$; $a_{mr} = 0$ otherwise. |
| $d_{rs}$ | $d_{rs} = 1$ if (sub-) transportation tasks $r$ and $s$ are decomposed from the same transportation task; $d_{rs} = 0$ otherwise. |
| $L_r^{\text{Route}}$ | Mileage of transportation task $r$. |
| $L_m^{\text{Interval}}$ | Mileage cycle of the second-level maintenance for train $m$, i.e. the mileage interval between two adjacent second-level maintenances. |
| $T_m^{\text{Begin}}$ | The starting planning time of train $m$ within the planning period. For a train that newly puts into operation or finishes a high-level maintenance, $T_m^{\text{Begin}}$ equals the day it puts into operation or the next day it finishes the high-level maintenance; for other trains, $T_m^{\text{Begin}}$ equals to the starting day of the planning period. |
| $T_m^{\text{End}}$ | The ending planning time of train $m$ within the planning period. For |
|   | a train that will reach its high-level maintenance cycle, $T_m^{\text{End}}$ equals the previous day it sends to the high-level maintenance; for other trains, $T_m^{\text{End}}$ equals to the ending day of the planning period. |
| $\Delta_p^{\text{Maintenance}}$ | Duration time of maintenance packet $r$. According to current practice, $\Delta_p^{\text{Maintenance}}$ is set as three days regardless of types of maintenance packets. |
| $C^{\text{Maintenance}}$ | Depot maintenance capacity, i.e. the maximal number of trains a depot can hold at the same day. |
| **Decision Variables** | |
| $x_m^r(t)$ | Binary integer, if train $m$ executes transportation task $r$ on day $t$, $x_m^r(t)=1$; otherwise, $x_m^r(t)=0$. |
| $y_m^p(t)$ | Binary integer, if train $m$ is sent to perform maintenance packet $r$ on day $t$, $y_m^p(t)=1$; otherwise, $y_m^p(t)=0$. |

## 4.3 Model formulation

Given the above-mentioned assumptions and notations, the high-speed train assignment and second-level maintenance planning problem can be formulated as a binary non-linear programing model as follows:

*Objective function*:

$$\min \sum_{m\in E} \sum_{t\in(T_m^{\text{Begin}},T_m^{\text{End}})} \sum_{p\in P} [L_m^{\text{Interval}} - l_m(t-1)] b_{mp} y_m^p(t) \tag{1}$$

*Subject to*:

$$\sum_{r\in R} a_{mr} x_m^r(t) + f_m(Y_m,t) = 1 \quad \forall m\in E,\ t\in(T_m^{\text{Begin}},T_m^{\text{End}}) \tag{2}$$

$$\sum_{m\in E} a_{mr} x_m^r(t) = 1 \quad \forall r\in R,\ t\in(0,T) \tag{3}$$

$$\sum_{m\in E} f_m(Y_m,t) \leq C^{\text{Maintenance}} \quad \forall t\in(0,T) \tag{4}$$

$$l_m(t) \leq L_m^{\text{Interval}} \quad \forall m\in E,\ t\in(T_m^{\text{Begin}},T_m^{\text{End}}) \tag{5}$$

$$\sum_{t\in(T_m^{\text{Begin}},T_m^{\text{End}})} \sum_{p\in P} b_{mp} y_m^p(t) \leq 1 \quad \forall m\in E \tag{6}$$

$$x_m^r(t) = x_m^{r+1}(t+1) \quad \forall m\in E,\ t\in(T_m^{\text{Begin}},T_m^{\text{End}}-1),\ d_{r,r+1}=1, r,r+1\in R \tag{7}$$

$$x_m^r(t),\ y_m^p(t) \in \{0,1\} \quad \forall m\in E,\ t\in(T_m^{\text{Begin}},T_m^{\text{End}}),\ r\in R,\ p\in P \tag{8}$$

where $l_m(t)$ is the cumulative mileage (from last second-level maintenance) for train

$m$ day $t$, and it is updated with the following equation:

$$l_m(t) = [l_m(t-1) + \sum_{r \in R} \delta_m^r L_r^{\text{Route}} x_m^r(t)][1 - \sum_{p \in P} y_m^p(t)] \tag{9}$$

While $f_m(Y_m, t)$ is the state function indicating whether a train $m$ is in the operation state ($f_m(Y_m, t) = 1$) or in the maintenance state ($f_m(Y_m, t) = 0$) on day $t$. The function can be expressed by:

$$f_m(Y_m, t) = \sum_{\tau \in (T_m^{\text{Begin}}, T_m^{\text{Begin}})} \sum_{p \in P} b_{mp} y_m^p(\tau) \phi(t \in [\tau, \tau + \Delta_p^{\text{Maintenance}}]) \tag{10}$$

where $Y_m$ represents the solution vector associated with train $m$, and

$$\phi(t \in S) = \begin{cases} 1 & \text{If } t \in S \\ 0 & \text{Otherwise} \end{cases} \tag{11}$$

Formula (1) is the objective function, which minimizes the mileage losses for all high-speed trains. Here the mileage loss means the theoretical mileage cycle of the second-level maintenance minus the cumulative operating mileage when the train is sent to the second-level maintenance. Constraint (2) ensures that a train either executes a transportation task or performs the second-level maintenance within the planning horizon. Constraint (3) guarantees that each transportation task is assigned exactly to one train. Constraint (4) confirms that the total number of trains performing the second-level maintenance on a particular day does not exceed the depot maintenance capacity. While Constraint (5) ensures that the cumulative mileage does not exceed the maintenance cycle. Constraint (6) guarantees that within the planning horizon a high-speed train carries out the second-level maintenance at most once (Assumption 1). Constraint (7) states that if two (sub-) transportation tasks are decomposed from a same transportation task, they must be connected by the same train. Finally, the variable domains are specified by Constraint (8).

### 4.4 Model discussions

This sub-section discusses two special scenarios: the multi-day (more than three days) routes and high-level maintenance considerations. Mathematical formulations are presented to model these scenarios.

### 4.4.1 Extension to multi-day routes

If a route lasts for more than three days, the model need to be extended with respect to Constraint (7). Let $r, r+1, r+2, \cdots, r+n-1$ denote the sequence of sub-routes decomposed from the $n$-day ($n \geq 3$) father route. If $x_m^r(t) = 1$ for $m \in E$, $t \in (T_m^{\text{Begin}}, T_m^{\text{End}} - n + 1)$, Constraint (7) should be modified as follows:

$$x_m^{r+1}(t+1) = x_m^{r+2}(t+2) = \cdots = x_m^{r+n-1}(t+n-1) = 1 \tag{12}$$

Though currently there exits no route that lasts for more than three days in the China high-speed system, we believe in the future a longer route may be employed as the high-speed network is able to connect two cities that locate significantly far apart from each other.

**4.4.2 High-level maintenance considerations**

Compared with the short-term tactical maintenance plan discussed in this paper, a high-level maintenance plan (HMP) is of long planning horizon nature and is designed at a strategical level, which aims to balance the maintenance work and avoid concentrated maintenance [21]. Once a HMP is established, the short-term high-speed train assignment and maintenance schedule should be designed with considerations of the HMP to ensure the HMP can be successfully implemented. For example, the HMP indicates a daily average running mileage for the high-speed trains based on which the trains can meet or approach the maintenance cycle on scheduled maintenance days as much as possible. To this end, railway practitioners always set a reasonable range of running mileage during the planning horizon for the trains. This rule can be mathematical described by the following formula:

$$L_m^{\min} \leq \sum_{t \in (T_m^{\text{Begin}}, T_m^{\text{End}})} \sum_{r \in R} \delta_m^r L_r^{\text{Route}} x_m^r(t) \leq L_m^{\max} \quad \forall m \in E \tag{13}$$

where $L_m^{\min}$ and $L_m^{\max}$ denote the lower and upper bound for the running mileage range of train $m$, respectively.

## 5. Conclusions

Assigning rolling stock into reasonable train routes while taking into account maintenance constraints is a fundamental problem faced by a passenger railway operator. This paper considers the simultaneous high-speed train assignment and maintenance scheduling problem. For the maintenance scheduling module, our focus is the second-level maintenance, which is carried out once a month on average. We propose a binary non-linear programming model to mathematically describe the problem. The objective aims to minimize the mileage losses for all high-speed trains and the constraints cover various operational requirements and capacity restrictions. To model the difficult operational requirements, a novel cumulative mileage update function is developed; meanwhile, to describe the depot maintenance capacity restriction, we employ a state function that is able to identify whether a train is in the maintenance state or in the operation state.

Our future research direction is to develop effective and efficient solution algorithms to solve the proposed mathematical model and apply the solution approach into real-world problem instances.